\documentclass[12pt]{amsart}
\usepackage{amsfonts,color,latexsym,amssymb,amsmath,epsfig,rotating}

\newtheorem{sz}{The Schwartz-Zippel Lemma}

\newtheorem{mink}{Minkowski's Theorem}

\newcommand{\res}{{\mathrm{Res}}} 
\newtheorem{dfn}{Definition}[section]

\newtheorem{rem}[dfn]{Remark} 
\newtheorem{prop}[dfn]{Proposition}

\newtheorem{thm}[dfn]{Theorem} 
\newtheorem{lemma}[dfn]{Lemma}
\newtheorem{cor}[dfn]{Corollary}

\newtheorem{ex}[dfn]{Example}

\renewcommand{\qed}{$\blacksquare$}

\newlength{\jmr}
\newcommand{\thth}{^{\text{\underline{th}}}}

\newcommand{\np}{{\mathbf{NP}}}

\newcommand{\bpp}{{\mathbf{BPP}}}
\newcommand{\zpp}{{\mathbf{ZPP}}}
\newcommand{\pp}{\mathbf{P}}
\newcommand{\ppoly}{\mathbf{P/poly}}

\newcommand{\nexp}{\mathbf{NEXP}}

\newcommand{\F}{\mathbb{F}}
\newcommand{\bF}{\overline{\F}}
\newcommand{\Q}{\mathbb{Q}}
\newcommand{\R}{\mathbb{R}}

\newcommand{\N}{\mathbb{N}}
\newcommand{\Z}{\mathbb{Z}}

\newcommand{\sat}{{\tt 3CNFSAT}} 
\newcommand{\cA}{{\mathcal{A}}}

\newcommand{\cS}{{\mathcal{S}}}

\newcommand{\cL}{{\mathcal{L}}}

\newcommand{\dia}{$\diamond$}

\newcommand{\size}{\mathrm{size}}

\author{Jingguo Bi}
\address{School of Mathematics, Shandong University, 
Jinan, 250100, P.R.\ China. } 
\email{jguobi@mail.sdu.edu.cn}
\author{Qi Cheng}
\address{School of Computer Science, University of Oklahoma, Norman, 
Oklahoma \ 73019, USA.}
\email{qcheng@cs.ou.edu} 
\author{J.\ Maurice Rojas}
\address{Department of Mathematics,
Texas A\&M University TAMU 3368,
College Station, Texas \ 77843-3368, USA.}  
\email{rojas@math.tamu.edu} 
\thanks{\mbox{}\\  
{\sc School of Mathematics, Shandong University, 
Jinan, 250100, P.R.\ China. }\\ 
{\em E-mail address:} {\tt jguobi@mail.sdu.edu.cn}\\  
J.B.\ was partially supported by NSF of China Projects 
(No.\ 61133013 and No.\ 60931160442).\\  
\mbox{}\\
{\sc School of Computer Science, University of Oklahoma, Norman, 
OK \ 73019, USA.}\\ 
{\sc E-mail:} {\tt qcheng@cs.ou.edu}\\ 
Q.C.\ was partially supported by NSF under grants CCF-0830522 and 
CCF-0830524.\\   
\mbox{}\\
{\sc Department of Mathematics, Texas A\&M University TAMU 3368,
College Station, Texas \ 77843-3368, USA.}\\ 
{\em E-mail address:} {\tt rojas@math.tamu.edu}\\  
J.M.R.\ was partially supported by NSF MCS grant DMS-0915245,
DOE ASCR grant DE-SC0002505, and Sandia National Laboratories.} 
\title[Sub-Linear Root Detection Over Finite Fields]{\mbox{}\\  
\vspace{-1in}Sub-Linear Root Detection, and New Hardness Results, for 
Sparse Polynomials Over Finite Fields}  

\begin{document}


\begin{abstract}   
We present a deterministic $2^{O(t)}q^{\frac{t-2}{t-1}+o(1)}$ algorithm  
to decide whether a univariate polynomial $f$, with exactly $t$ monomial 
terms and degree $<\!q$, has a root in $\F_q$. A corollary of our 
method --- the first with complexity {\em sub-linear} in $q$ when $t$ is 
fixed --- is that the nonzero roots in $\F_q$
can be partitioned into no more than $2\sqrt{t-1}(q-1)^{\frac{t-2}{t-1}}$ 
cosets of two proper subgroups $S_1 \subseteq S_2$ of $\F^*_q$.  
Another corollary is the first 
deterministic sub-linear algorithm for detecting common degree one factors of 
$k$-tuples of $t$-nomials in $\F_q[x]$ 
when $k$ and $t$ are fixed. 

When $t$ is not fixed we show that each of the following 
problems is $\np$-hard with respect to $\bpp$-reductions, even 
when $p$ is prime: 
\begin{itemize}
\item[$\bullet$]{detecting roots in $\F_p$ for $f$}  
\item[$\bullet$]{deciding whether the square of a degree one polynomial in 
$\F_p[x]$ divides $f$}
\item[$\bullet$]{deciding whether the square of a degree one polynomial in
$\bF_p[x]$ divides $f$}
\item[$\bullet$]{deciding whether the gcd of two $t$-nomials in
$\F_p[x]$ has positive degree} 
\end{itemize} 
Finally, we prove that if the 
complexity of root detection is sub-linear (in a refined sense),  
relative to the {\em straight-line program encoding}, then 
$\nexp\!\not\subseteq\!\ppoly$. 
\end{abstract} 

\maketitle 
\titlepage 

\addtocounter{page}{1} 

\section{Introduction} 
The solvability of univariate sparse polynomials is a fundamental problem
in computer algebra, and an important precursor to deep questions in 
polynomial system solving and circuit complexity. Cucker, Koiran, and Smale 
\cite{cks} found a polynomial-time algorithm to \mbox{find} all integer roots 
of a univariate polynomial $f$ in $\Z[x]$ with exactly $t$ terms, i.e., a 
{\em univariate $t$-nomial}. Shortly 
afterward, H.\ W.\ Lenstra, Jr.\ \cite{lenstrafac} gave a polynomial-time 
\mbox{algorithm} to compute all factors of fixed degree over an algebraic 
extension of $\Q$ of fixed degree (and thereby all rational roots).
Independently, Kaltofen and Koiran \cite{kaltofenkoiran} and Avendano, 
Krick, and Sombra \cite{avendanokricksombra} extended this to finding 
bounded-degree factors of sparse polynomials in $\mathbb{Q}[x,y]$ in 
polynomial-time. Unlike the famous LLL factoring algorithm \cite{lll}, the 
complexity for the algorithms from 
\cite{cks,lenstrafac,kaltofenkoiran,avendanokricksombra} was 
relative\linebreak 
\scalebox{.93}[1]{to the {\em sparse encoding} (cf.\ Definition 
\ref{dfn:sparse} of Section \ref{sec:back} below) and thus polynomial in 
$t+\log \deg f$.}  

Changing the ground field dramatically changes the complexity. 
For instance, while polynomial-time algorithms are now known for detecting real 
roots for trinomials in $\Z[x]$ \cite{rojasye,brs}, no polynomial-time 
algorithm is known for tetranomials \cite{bhpr}. Also, 
detecting $p$-adic rational 
roots for trinomials in $\Z[x]$ was only recently shown to lie in $\np$ (for 
fixed $p$), as was $\np$-hardness with respect to $\zpp$-reductions for 
$t$-nomials when neither $t$ nor $p$ are fixed \cite[Thm.\ 1.4 \& Cor.\ 1.5]
{airr}.    

Here, we focus on the complexity of detecting solutions of univariate 
$t$-nomials over finite fields. 

\subsection{Main Results and Related Work}   
While deciding the existence of a $d\thth$ 
root of an element of the $q$-element field $\F_q$ is doable in time 
polynomial in $\log(d)+\log q$ (see, e.g., \cite[Thms.\ 5.6.2 \& 5.7.2, 
pg.\ 109]{bs}), detecting roots for a {\em trinomial} equation 
$a+bx^{d_0}+cx^d\!=\!0$ with $d\!>\!d_0\!>\!0$ within time sub-linear in 
$d$ and $q$ is already a mystery. 
(Erich Kaltofen and David A.\ Cox independently asked about such 
polynomial-time algorithms around 2003 \cite{kaltofen,cox}.)  
We make progress on a natural extension of this question. In what 
follows, we use $|S|$ for the cardinality of a set $S$.  
\begin{thm} 
\label{thm:meat} 
Given any univariate $t$-nomial  
$$f(x) := c_1 + c_2 x^{a_2} + c_3 x^{a_3} + \cdots + c_t x^{a_t}\in\F_q[x]$$  
with degree $<\!q$, we can decide, within $4^{t}(t\log q)^{O(1)}
+t^{\frac{1}{2}+o(1)}q^{\frac{t-2}{t-1}+o(1)}$
deterministic bit\linebreak operations, whether $f$ has a root in $\F_q$. 
Moreover, letting $\delta\!:=\!\gcd(q-1,a_2,\ldots,a_t)$ and 
$\eta\!:=\!\sqrt{t-1}\left(\frac{q-1}{\delta}\right)^{\frac{t-2}{t-1}}$, 
the entire set of nonzero roots of $f$ in $\F_q$ is a union of at most 
$2\eta$ cosets of two proper subgroups $S_1 \subseteq S_2$ of $\F_q^*$, 
where $|S_1|\!=\!\delta$ and $\frac{\delta^{\frac{t-2}{t-1}}(q-1)^{\frac{1}
{t-1}}}{\sqrt{t-1}}\!\leq\!|S_2|\!\leq\!\frac{q-1}
{2}$. In particular, the number of nonzero roots of $f$ is no more than  
$\max\left\{2\delta\eta,\frac{\eta-1}{\eta}\cdot \frac{q-1}{2}\right\}$. 
\end{thm} 

\noindent 
The degree assumption is natural since $x^q\!=\!x$ in 
$\F_q[x]$. Note also that deciding whether an $f$ as above has a root in $\F_q$ 
via brute-force search takes $q^{1+o(1)}$ bit operations, assuming $t$ is 
fixed. 

Our first main result thus includes a finite 
field analogue of {\em Descartes' Rule} \cite{descartes}. (The latter result 
implies an upper bound of $2t-1$ for the number of real roots of a real 
univariate $t$-nomial.) More to the point, Theorem \ref{thm:meat} provides new 
structural and algorithmic information, complementing an earlier 
finite field analogue of Descartes' Rule \cite[Lemma 7]{cfklls}. Theorem 
\ref{thm:meat} can also be thought of as a refined, positive characteristic 
analogue of results of Tao and Meshulam \cite{tao,meshulam} bounding the 
number of complex roots of unity at which a sparse polynomial can vanish 
(a.k.a.\ uncertainty inequalities over finite groups). 

Note that if we pick $a_2,\ldots,a_t$ uniformly randomly in $\{-M,\ldots,M\}$ 
then, as $M\!\longrightarrow \infty$, the probability that 
$\gcd(a_2,\cdots,a_t)\!=\!1$ approaches $1/\zeta(t-1)$ (see, e.g., 
\cite{chris}). The latter quantity increases from 
$\frac{6}{\pi^2}\!\approx\!0.6079$ to $1$ as $t$ goes from $3$ to $\infty$. 
Our theorem thus implies that, with ``high'' probability, the rational 
roots of a sparse polynomial over a finite field can be divided into
two components: one component consisting of few isolated roots, and the other 
component consisting of few cosets of a (potentially large) subgroup
of $\F_q^*$. Put another way, 
if the number of the rational 
roots of a sparse polynomial is close to its degree, 
then the set of the roots must exhibit a strong multiplicative structure.

Since detecting roots over $\F_q$ is the same as detecting linear factors of 
polynomials in $\F_q[x]$, it is natural to ask about the complexity of 
factoring sparse polynomials over $\F_q[x]$. The asymptotically fastest 
randomized algorithm for factoring arbitrary $f\!\in\!\F_q[x]$ of 
degree $d$ uses $O(d^{1.5}+d^{1+o(1)}\log q)$ arithmetic operations 
in $\F_q$ \cite{kedlaya}, but no complexity bound polynomial in 
$t+\log(d)+\log q$ is known. 
(See \cite{berlekamp,cantorzassenhaus,kaltofenshoup,umans} for some important 
milestones, and \cite{gp,kaltofen,gat} for an extensive survey on factoring.) 
However, to detect roots in $\F_q$, we don't need the full power of factoring: 
we need only decide whether $\gcd(x^q-x,f(x))$ has positive degree. 
Indeed, a consequence of our first main result is a speed-up for a 
variant of the latter decision problem. 
\begin{cor}  
\label{cor:kby1}
Given any univariate $t$-nomials $f_1,\ldots,f_k\!\in\!\F_q[x]$,  
we can decide 
if $f_1,\ldots,f_k$ have a common degree one  
factor in $\F_q[x]$ via a deterministic algorithm  
with complexity 
$4^{kt-k}(kt\log q)^{O(1)}
+k\left(k\sqrt{t}\right)^{1+o(1)}q^{\frac{kt-k-1}{kt-k}+o(1)}$.  
\end{cor} 

\noindent 
Corollary \ref{cor:kby1} appears to give the first sub-linear algorithm 
for detecting roots of $k$-tuples of univariate $t$-nomials for 
$k$ and $t$ fixed. 
\begin{rem} 
It is important to note that the $k\!=\!2$ case is {\em not} the same 
as deciding whether the gcd of two {\em general} polynomials 
has positive degree: the latter problem is the same as detecting common 
factors of {\em arbitrary} degree, or degree one factors over an extension 
field. Finding an algorithm for the latter problem with 
complexity sub-linear in $q$ is already an open problem for $k\!=\!2$ and 
$t\!\geq\!3$: see \cite{emirispan}, and Theorem  
\ref{thm:adisc} and Remark \ref{rem:res} below. \dia 
\end{rem} 

One reason why it is challenging to attain complexity sub-linear 
in $q$ is that detecting roots in $\F_q$ for $t$-nomials 
is $\np$-hard when $t$ is not fixed, even restricting to one 
variable and prime $q$.  
\begin{thm}
\label{thm:fphard} 
Suppose that, for any input $(f,p)$ with $p$ a prime and $f\!\in\!\F_p[x]$ 
a $t$-nomial of degree $<p$, one could decide whether $f$ has a root in 
$\F_p$ within $\bpp$, using $t+\log p$ as the underlying input size. Then 
$\np\!\subseteq\!\bpp$. 
\end{thm} 

\noindent 
The least $n$ making root detection in $\F^n_p$ be $\np$-hard for 
polynomials in $\F_p[x_1,\ldots,x_n]$ (for $p$ prime, and relative to the 
sparse encoding) appears to have been unknown. Theorem \ref{thm:fphard} thus 
comes close to settling this problem.  
Theorem \ref{thm:fphard} also complements an earlier result of 
Kipnis and Shamir proving $\np$-hardness 
for detecting roots of univariate sparse polynomials over fields of the form 
$\F_{2^\ell}$ \cite{kipnisshamir}. Furthermore, Theorem \ref{thm:fphard} 
improves another recent $\np$-hardness result where the underlying 
input size was instead the (smaller) straight-line program complexity 
\cite{chenghillwan}. 

Let $\bF_q$ denote the algebraic closure of $\F_q$. 
A consequence of our last complexity lower bound is the hardness of detecting 
{\em degenerate} roots over $\F_p$ and $\bF_q$:  
\begin{thm} 
\label{thm:adisc} 
Consider the following two problems, each with input $(f,p)$ where $p$ is a 
prime and $f\!\in\!\F_p[x]$ is a $t$-nomial of degree $<p$. 
\begin{enumerate} 
\item{Decide whether $f$ is divisible by the square of a degree one polynomial 
in $\F_p[x]$.} 
\item{Decide whether $f$ is divisible by the square of a degree one polynomial
in $\bF_p[x]$.}
\end{enumerate} 
Then, using $t+\log p$ as the underlying input size, each of these problems is 
$\np$-hard with respect to $\bpp$-reductions. 
\end{thm}  

\noindent 
The $\np$-hardness of both problems had been previously unknown. 
Theorem \ref{thm:adisc} thus improves \cite[Cor.\ 2]{squarefree} where  
$\np$-hardness (with respect to $\bpp$-reductions)  
was proved for the {\em harder} variant of Problem (2) where one expands the 
allowable inputs to polynomials in $\bF_p[x]$. 
\begin{rem} 
Note that detecting a degenerate root for $f$ is the same as 
detecting a common degree one factor of $f$ and $\frac{\partial f}{\partial 
x}$, at least when $\deg f$ is less than the characteristic of the field.  
So an immediate consequence of Theorem \ref{thm:adisc} is that 
detecting common degree one factors in $\F_p[x]$ (resp.\ $\bF_p[x]$) 
for pairs of polynomials in $\F_p[x]$ is $\np$-hard with 
respect to $\bpp$-reductions. We thus also strengthen earlier work 
proving similar complexity lower bounds for detecting common degree one  
factors in $\F_q[x]$ (resp.\ $\bF_q[x]$) \cite[Thm.\ 4.11]{von}. \dia 
\end{rem}  
\begin{rem} 
\label{rem:res} 
It should be noted that Problem (2) is equivalent to deciding the 
vanishing of univariate {\em $\cA$-discriminants} (see \cite[Ch.\ 12, pp.\ 
403--408]{gkz94} and Definitions \ref{dfn:syl} and \ref{dfn:adisc} of 
Section \ref{sub:backhard} below). While Lemma \ref{lemma:tri} of Appendix A  
tells us that the trinomial case of Problem (2) can be 
done in $\pp$, we are unaware of any other speed-ups for fixed $t$. 
In particular, it follows immediately from Theorem \ref{thm:adisc} that 
deciding the vanishing of univariate {\em resultants} (see, e.g., 
\cite[Ch.\ 12, Sec.\ 1, pp.\ 397--402]{gkz94} and Definition \ref{dfn:syl} 
of Section \ref{sub:backhard} below), 
over $\F_p[x]$, is also $\np$-hard with respect to $\bpp$-reductions. \dia 
\end{rem}

Our final result is a complexity separation depending on a weak 
tractability assumption for detecting roots of univariate polynomials given 
as {\em straight-line programs (SLPs)}. 
\begin{thm}  
\label{thm:slp} 
Suppose that, given any straight-line program of size $L$ representing a 
polynomial $f\!\in\!\F_{2^\ell}[x]$, we could decide if $f$ has a root 
in $\F_{2^\ell}$ within time $L^{O(1)}2^{\ell-\omega(\ell)}$. Then 
$\nexp\!\not\subseteq\!\ppoly$. 
\end{thm} 

\noindent  
One should recall that $\nexp\!\subseteq\!\ppoly \Longleftrightarrow 
\nexp\!=\!\mathbf{MA}$ \cite{ikw}. So the conditional assertion 
of our last theorem indeed implies a new separation of complexity classes. 
It may actually be the case that there is no algorithm for detecting roots 
in $\F_{2^\ell}$ better than brute-force search. Such a result would be in 
line with the Exponential Time Hypothesis \cite{impagliazzo} 
and the widely-held belief in the 
cryptographic community that the 
only way to break a well-designed block cipher is by exhaustive search.  

\subsection{Highlights of Main Techniques} 
The key new advance needed to attain our speed-ups is a method, 
based on the {\em Shortest Vector Problem (SVP)} for a lattice basis 
(see \cite{mic} and Section \ref{sub:backspeed}),   
to lower the degree of any sparse polynomial in $\F_q[x]$ to a  
power of $q$ strictly less than $1$ while still preserving solvability 
over $\F_q$. 
\begin{lemma} 
\label{lemma:expo}
Given integers $a_1,\cdots, a_t,N$ satisfying 
$0\!<\!a_1\!<\cdots<\!a_t\!<\!N$  
and\linebreak $\gcd(N,a_1,\cdots, a_t) = 1$, one can find, within  
$4^{t}(t\log N)^{O(1)}$ bit operations, 
an integer $e$ with the following property for all $i\!\in\!\{1,\ldots,t\}$: 
if $m_i\!\in\!\left\{-\lfloor N/2\rfloor,\ldots,\lceil N/2 
\rceil\right\}$ is the unique integer congruent to $ea_i$ mod $N$ then
$|m_i|\!\leq\!\sqrt{t} N^{1-t^{-1}}$.  
\end{lemma}

\noindent 
We prove this lemma in Section \ref{sub:backspeed}, and show 
how the lemma can be applied to the 
exponents of a sparse polynomial to yield Theorem \ref{thm:meat}
in Section \ref{sub:meat}.  
Corollary \ref{cor:kby1} is proved in Section \ref{sub:kby1}. 
\begin{ex} 
\label{ex:tetra} 
Consider any polynomial of the form\\ 
\mbox{}\hfill  
$f(x)\!=\!c_1+c_2x+c_3x^{2^{200}+26}+c_4x^{2^{200}+27}\!\in\!\F_q[x]$\hfill 
\mbox{}\\ 
where\\  
\mbox{}\hfill 
$q\!:=\!6(2^{200}+26)+1\!=\!
9641628265553941653251772554046975615133217962696757011808413
$\hfill \mbox{}\\  
(which is a $61$-digit prime) and 
$c_1c_4\!\neq\!0$. Considering the lattice generated by the vectors 
$(1,2^{200}+26,2^{200}+27),(q-1,0,0),(0,q-1,0),(0,0,q-1)$, it is not hard to 
see that $(6,0,6)$ is a minimal length vector in this lattice. Moreover, 
$6\cdot 1\!\equiv\!6$, 
$6(2^{200}+26)\!\equiv\!0$, $6(2^{200}+27)\!\equiv\!6$ mod $q-1$. 
Letting $\sigma$ be any generator of $\F^*_q$ it is clear that any 
$x\!\in\!\F^*_q$ can be written as $x\!=\!\sigma^i z$ for some $i\!\in\!
\{0,\ldots,5\}$ and $z\!\in\!\F^*_q$ satisfying $z^{\frac{q-1}{6}}\!=\!1$. 
So then, we see that 
solving $f(x)\!=\!0$ is equivalent to 
finding an $i\!\in\!\{0,\ldots,5\}$  
and a $z\!\in\!\F^*_q$ with\\ 
\mbox{}\hfill $\left(c_1+c_3\sigma^{(2^{200}+26)i}\right)
+\left(c_2\sigma^i+c_4\sigma^{(2^{200}+27)i}\right)z^6 \ = \ 
z^{\frac{q-1}{6}}-1 \ = \ 0$. \dia \hfill \mbox{} 
\end{ex} 

Recall that any Boolean expression of one of the following forms:\\
\mbox{}\hfill $(\diamondsuit)$ $y_i\vee y_j \vee y_k$, \
$\neg y_i\vee y_j \vee y_k$, \
$\neg y_i\vee \neg y_j \vee y_k$,  \
$\neg y_i\vee \neg y_j \vee \neg y_k$,
with $i,j,k\!\in\![3n]$,\hfill \mbox{} \\
is a $\sat$ {\em clause}.
A {\em satisfying assigment} for an
arbitrary Boolean formula $B(y_1,\ldots,y_n)$ is
an assigment of values from $\{0,1\}$ to the variables
$y_1,\ldots,y_{n}$ which makes the equality $B(y_1,\ldots,y_n)\!=\!1$
true.\footnote{We respectively identify $0$ and $1$ with ``False'' and
``True''.}
 
A key construction behind the proofs of Theorems \ref{thm:fphard} and 
\ref{thm:adisc} in Section \ref{sec:primes} is a highly structured 
randomized reduction from $\sat$ to detecting roots of 
univariate polynomial systems over finite fields. In particular, 
the finite fields arising in this reduction have cardinality coming from a 
very particular family of prime numbers. (See Definition \ref{dfn:sparse} 
from Section \ref{sec:back} for our definition of input size.) 
\begin{thm}
\label{thm:newplai}
Given any $\sat$ instance $B(y_1,\ldots,y_n)$ in $n\!\geq\!4$ variables with
$k$ clauses, there is
a (Las Vegas) randomized polynomial-time algorithm that produces
positive integers $c,p_1,\ldots,p_n$ and a
$k$-tuple of polynomials $(f_1,\ldots,f_k)\!\in\!\Z[x]$ with the
following properties:
\begin{enumerate}
\item{$c\!\geq\!11$ and $\log(cp_1\cdots p_n)\!=\!n^{O(1)}$.}
\item{$p_1,\ldots,p_n$ is an increasing sequence of primes
and $p\!:=\!1+cp_1\cdots p_n$ is prime.}
\item{For all $i$, $f_i$ is monic, $f_i(0)\!\neq\!0$, 
$\deg f_i\!<\!p_1\cdots p_n$, and $\size(f_i)\!=\!n^{O(1)}$.}
\item{For all $i$, the mod $p$ reduction of $f_i$ has exactly
$\deg f_i$ distinct roots in $\F_p$.}
\item{$B$ has a satisfying assignment if and only if
the mod $p$ reduction of $(f_1,\ldots,f_k)$ has a root in $\F_p$. \qed }
\end{enumerate}
\end{thm}

\noindent 
Theorem \ref{thm:newplai} is based on an earlier reduction of  
Plaisted involving complex roots of unity 
\cite[Sec.\ 3, pp.\ 127--129]{plaisted} and was refined into the form
below in \cite[Secs.\ 6.2--6.3]{airr}.\footnote{\cite{airr} in fact contains
a version of Theorem \ref{thm:newplai} with $c\!\geq\!2$, but $c\!\geq\!11$
can be attained by a trivial modification of the proof there.}

We now review some additional background necessary for our proofs. 

\section{Background} 
\label{sec:back}  
Our main notion of input size essentially reduces to
how long it takes to write down monomial term expansions, a.k.a. the
{\em sparse encoding}.
\begin{dfn}
\label{dfn:sparse}
For any polynomial $f\!\in\!\Z[x_1,\ldots,x_n]$ written
$f(x)\!=\!\sum^t_{i=1} c_i x^{a_{1,i}}_1
\cdots x^{a_{n,i}}_n$, we define $\size(f)\!:=\!\sum^t_{i=1}
\log_2\left[(2+|c_i|)(2+|a_{1,i}|)\cdots
(2+|a_{n,i}|)\right]$. Also, when $F\!:=\!(f_1,\ldots,f_k)$, we define
$\size(F)\!:=\!\sum^k_{i=1}\size(f_i)$. \dia
\end{dfn}

\noindent
The definition above is also sometimes known as the {\em sparse} size
of a polynomial. Note that $\size(c)\!=\!O(\log |c|)$ for any integer
$c$.  

A useful fact, easily obtainable from the famous {\em Schwartz-Zippel Lemma} 
is that systems of univariate polynomial equations can, at the expense of some 
randomization, be reduced to {\em pairs} of univariate equations. (See 
Appendix A for the proof and \cite{giustiheintz} for a multivariate version.) 
\begin{lemma}
\label{lemma:lin}
Given any prime power $q$ and $f_1,\ldots,f_k\!\in\!\F_q[x]$,
let $Z(f_1,\ldots,f_k)$ denote the set of solutions of
$f_1\!=\cdots=\!f_k\!=\!0$ in $\bF_q$. Also let $d\!:=\!\max_i \deg f_i$.
Then at least a fraction of $1-\frac{d}{q}$ of the
$(u_2,\ldots,u_k)\!\in\!\F^{k-1}_q$ satisfy
$Z(f_1,\ldots,f_k)\!=\!Z(f_1,u_2f_2+\cdots+u_kf_k)$.
\end{lemma}  
\begin{rem} 
\label{rem:lingood} 
For this lemma to yield a high-probability reduction from $k\times 1$ 
systems to $2\times 1$ systems, we will of course need to assume 
that $d$ is a small constant fraction of $q$. This will indeed be the 
case in our upcoming applications since we will be combining 
the lemma\linebreak 
\scalebox{.98}[1]{with Theorem \ref{thm:newplai}, and Assertions (1)--(3) of 
the theorem force $d\!<\!\frac{p}{11}$ (with $q\!=\!p$ a prime). \dia } 
\end{rem} 

Let us now observe the following complexity bound for root detection for
(not necessarily sparse) polynomials over finite fields.
\begin{prop}
\label{prop:fac}
Given any polynomial $f\!\in\!\F_q[x]$ of degree $d$ and $N|(q-1)$, 
we can decide within $d^{1+o(1)} (\log q)^{2+o(1)}$ deterministic 
bit operations whether $f$ has a root in the order $N$ subgroup 
of $\F^*_q$. \qed
\end{prop}

\noindent 
Since detecting roots for $f$ as above is the same as deciding whether
$\gcd(x^N-1,f(x))$ has positive degree, the complexity bound
above can be attained as follows: compute $r(x)\!:=\!x^N$ mod $f(x)$ via
recursive squaring \cite[Thm.\ 5.4.1, pg.\ 103]{bs}, and then
compute\linebreak $\gcd(r(x)-1,f(x))$ in time $d^{1+o(1)}(\log q)^{1+o(1)}$ 
via the Knuth-Sch\"onhage algorithm \cite[Ch.\ 3]{bcs}.  

\subsection{Geometry of Numbers for Speed-Ups} 
\label{sub:backspeed} 
Recall that a {\em lattice} 
in $\mathbb{R}^{m}$ is the set $\cL(\mathbf{b_1},\ldots,\mathbf{b_d})=
\left\{\left. \sum\limits^d_{i = 1}x_i\mathbf{b_i}\; \right| 
\; x_i\in \mathbb{Z}\right\}$ 
of all integral combinations of $d$ linearly independent vectors
$\mathbf{b_1},\ldots,\mathbf{b_d}\in \mathbb{R}^{m}$. The integers $d$ and $m$
are respectively called the {\em rank} and {\em dimension} of the lattice.
The determinant $\det(\cL)$ of the lattice $\cL$ is
the volume of the $d$-dimensional parallelepiped spanned by the origin and the
vectors of any $\Z$-basis for $\cL$. Any lattice can be conveniently 
represented by a $d\times m$ 
matrix $\mathbf{B}$, where $\mathbf{b_1},\ldots,\mathbf{b_d}$ are the rows.  
The determinant of the lattice $\cL$
can then be computed as $\det(\cL(\mathbf{B}))\!=\!
\sqrt{\det(\mathbf{ B}\mathbf{B}^\top)}$. 

Let $\| \cdot \|$ denote the Euclidean norm on $\R^n$ 
for any $n$. Perhaps the most famous computational problem on lattices is the 
{\em (exact) Shortest Vector Problem (SVP)}: Given a basis of a lattice $\cL$, 
find a non-zero vector $\mathbf{u}\!\in\!\cL$, such that 
$\| \mathbf{v}\| \geq \| \mathbf{u}\|$ for any  vector 
$\mathbf{v}\in \cL\setminus {\mathbf{0}}$. The following is a well-known upper 
bound on the shortest vector length in lattice $\cL$. 
\begin{mink}  
Any lattice $\cL$ of rank $d$ contains a non-zero
vector $\mathbf{v}$ with\linebreak $\|\mathbf{v}\|\!\leq\!\sqrt{d}  
\det(L)^{1/d}$. \qed 
\end{mink}  

Given a lattice with rank $d$, the celebrated {\em LLL algorithm} \cite{lll} 
can find, in time polynomial in the bit-size of a given basis for $\cL$, a 
vector whose length is at most $2^{\frac{d}{2}}$ times the length of the 
shortest nonzero vector in $\cL$. An algorithm 
with arithmetic  
complexity $d^{O(1)}4^d$,  
proposed in \cite[Sec.\ 5]{mic} by Micciancio and 
Voulgaris, is currently the fastest deterministic algorithm for 
solving SVP. (See \cite{ngu} for a survey of other SVP algorithms.) 

Let us now prepare for our degree-lowering tricks. First, we 
construct the lattice $\cL$ spanned by the rows of matrix 
$\mathbf{B}$, where\\  
($\star\star$) \hfill 
$\mathbf{B} = \begin{bmatrix}
    a_1 & a_2 & \cdots & a_{t}\\
    N & 0 & \cdots & 0 \\
    0 & N & \cdots & 0 \\
    \vdots &\vdots & \ddots & 0 \\
    0 & 0 & \cdots & N \\
  \end{bmatrix}$ \hfill \mbox{}\\ 
Letting $\mathbf{v}\!:=\!(m_1,m_2,\cdots,m_{t})$ be the shortest vector of 
lattice $\cL$, there then clearly exists an integer $e$ such that
$ea_1\!\equiv\!m_1,\ldots,ea_{t} \equiv m_{t}$ mod $N$. (In fact, 
$e$ is merely the coefficient of $(a_1,\ldots,a_t)$ in the underlying 
linear combination defining $\mathbf{v}$.) Most importantly, 
the factorization of $\det(\cL)$ is rather restricted when the $a_i$ are 
relatively prime.  
\begin{lemma} 
\label{lemma:det}
If $\gcd(N,a_1,\ldots,a_t)\!=\!1$ then $\det(\cL)|N^{t-1}$. 
\end{lemma}  

\noindent
{\bf Proof:} Let $\cL_i$ denote the sublattice of $\cL$ generated 
by all rows of $\mathbf{B}$ save the $i\thth$ row. Clearly then, 
$\det(\cL)|\det(\cL_i)$ for all $i$. Moreover, we have 
$\det(\cL_1)\!=\!N^t$ and, via minor expansion from the 
$i\thth$ column of $B$, we have 
$\det(\cL_{i+1})\!=\!a_iN^{t-1}$ for all $i\!\in\!\{1,\ldots,t\}$. 
So $\det(\cL)$ divides $a_1N^{t-1},\ldots,a_tN^{t-1}$ and we are done. \qed 

\smallskip 
We are now ready to prove Lemma \ref{lemma:expo}.  

\smallskip 
\noindent 
{\bf Proof of Lemma \ref{lemma:expo}:}
From Lemma \ref{lemma:det} and Minkowski's theorem, there exists a shortest 
vector $\mathbf{v}$ of $\cL$ satisfying $\|\mathbf{v}\|\! 
\leq\!\sqrt{t}N^{1-t^{-1}}$. 
By invoking the exact SVP algorithm from \cite{mic} we can then find the 
shortest vector $\mathbf{v}$ in time $4^t (t\log N)^{O(1)}$. 
Let $\mathbf{v}\!:=\!(m_1,\ldots,m_{t})$. Clearly, by shortness, we 
may assume $|m_i|\!\leq\!N/2$ for all $i\!\in\!\{1,\ldots,t\}$. 
(Otherwise, we would be able to 
reduce $m_i$ in absolute value by subtracting a suitable row of the matrix 
$\mathbf{B}$ from $\mathbf{v}$.) Also, by construction, there is an $e$ such 
that $ea_i \equiv m_i$ mod $N$ for all $i\!\in\!\{1,\ldots,t\}$. \qed

\subsection{Resultants, $\cA$-discriminants, and Square-Freeness} 
\label{sub:backhard}  
Let us first recall the classical univariate resultant.
\begin{dfn}
\label{dfn:syl}
(See, e.g., \cite[Ch.\ 12, Sec.\ 1, pp.\ 397--402]{gkz94}.)
Suppose\linebreak $f(x)\!=\!a_0+\cdots+a_dx^{d}$ and
$g(x)\!=\!b_0+\cdots+b_{d'}x^{d'}$ are polynomials with
indeterminate coefficients. We define their {\em Sylvester matrix}
to be the $(d+d')\times (d+d')$ matrix

\noindent
\mbox{}\hfill \scalebox{1}[1]{$\cS_{(d,d')}(f,g)\!:=\!\begin{bmatrix}
a_0 & \cdots & a_d    & 0       & \cdots & 0 \\
   & \ddots &  & & \ddots &  \\
0   & \cdots & 0 & a_0 & \cdots & a_d \\
b_0 & \cdots & b_{d'}    & 0       & \cdots & 0 \\
  & \ddots &  &  & \ddots &  \\
0   & \cdots & 0 & b_0 & \cdots & b_{d'} 
\end{bmatrix}
\begin{matrix}
\\
\left. \rule{0cm}{.9cm}\right\}
d' \text{ rows}\\
\left. \rule{0cm}{.9cm}\right\}
d \text{ rows} \\
\\
\end{matrix}$}\hfill\mbox{}\\
and their {\em Sylvester resultant} to be
$\res_{(d,d')}(f,g)\!:=\!\det \cS_{(d,d')}(f,g)$. \dia
\end{dfn}
\begin{lemma}
\label{lemma:syl}
Following the notation of Definition \ref{dfn:syl},
assume $f,g\!\in\!K[x]$ for some field $K$, and that
$a_d$ and $b_{d'}$ are not both $0$. Then $f\!=\!g\!=\!0$ has a root in the
algebraic closure of $K$ if and only if $\res_{(d,d')}(f,g)\!=\!0$. More
precisely, we
have $\res_{(d,d')}(f,g)\!=\!a^{d'}_d \!\! \prod\limits_{f(\zeta)=0} g(\zeta)$
where the\linebreak

\vspace{-.5cm}
\noindent
product counts multiplicity. \qed
\end{lemma}
\noindent
The lemma is classical: see, e.g., \cite[Ch.\ 12, Sec.\ 1, pp.\ 
397--402]{gkz94}, \cite[pg.\ 9]{rs}, and \cite[Thm.\ 4.16, pg.\ 107]{bpr} 
for a more modern treatment.

We may now define a refinement of the classical {\em discriminant}. 
\begin{dfn}
\label{dfn:adisc} 
(See also \cite[Ch.\ 12, pp.\ 403--408]{gkz94}.)  
Let  $\cA\!:=\!\{a_1,\ldots,a_t\}\!\subset\!\N\cup \{0\}$
and $f(x)\!:=\!\sum^t_{i=1}c_ix^{a_i}$, where 
$0\!\leq\!a_1\!<\cdots<\!a_t$ and the $c_i$ are indeterminates. We then 
define the {\em $\cA$-discriminant} of $f$, $\Delta_\cA(f)$, to be\\
\mbox{}\hfill $\res_{(\bar{a}_t,\bar{a}_t-\bar{a}_2)}\left.
\left(\bar{f},\left.\frac{\partial\bar{f}}{\partial x}
\right/x^{\bar{a}_2-1}\right)\right/
c^{\bar{a}_t-\bar{a}_{t-1}}_t$, 
\hfill\mbox{}\\
where $\bar{a}_i\!:=\!(a_i-a_1)/g$ for all $i$, $\bar{f}(x)\!:=\!
\sum^t_{i=1}c_ix^{\bar{a}_i}$, and $g\!:=\!\gcd(a_2-a_1,\ldots,
a_t-a_1)$. \dia
\end{dfn}
\begin{rem}
Note that when $\cA\!=\!\{0,\ldots,d\}$ we have
$\Delta_\cA(f)\!=\!\res_{(d,d-1)}(f,f')/c_d$, i.e.,
for dense polynomials, the $\cA$-discriminant agrees with the
classical discriminant. \dia
\end{rem}

\begin{lemma}
\label{lemma:sq}
Suppose $p$ is any prime and $f,g\!\in\!\F_p[x]$ are relatively prime 
polynomials satisfying $f(0)g(0)\!\neq\!0$, $d\!:=\!\deg g\!\geq\!\deg f$, 
and $p\!>\!d$. Then the polynomial $f+ag$ is square-free for at least a 
fraction of $1-\frac{2d-1}{p}$ of the $a\!\in\!\F_p$.  
\end{lemma}
\begin{rem} 
\label{rem:sqgood} 
Just as for Lemma \ref{lemma:lin}, we will need to assume
that $d$ is a small constant fraction of $q$ for Lemma \ref{lemma:sq} 
to be useful. This will indeed be the case in our upcoming\linebreak 
applications since 
the setting will be the polynomials coming from Theorem \ref{thm:newplai}, 
and\linebreak Assertions (1)--(3) of the theorem force 
$2d-1\!<\!\frac{2}{11}p$ (with $q\!=\!p$ a prime). \dia 
\end{rem}

\noindent 
A stronger assertion, satisfied on a much smaller set of $a$, was 
observed earlier in the proof of Theorem 1 of \cite{squarefree}. For our 
purposes, easily finding an $a$ with $f+ag$ square-free will be crucial.   
We prove Lemma \ref{lemma:sq} in Appendix B. 

\section{Faster Root Detection: Proving Theorem \ref{thm:meat} and 
Corollary \ref{cor:kby1}}  
\label{sec:tnomial}

\subsection{Proving Theorem \ref{thm:meat}}   
\label{sub:meat} 
Before proving Theorem \ref{thm:meat}, let us first prove a result that 
will in fact enable sub-linear root detection in {\em arbitrary} subgroups 
of $\F^*_q$. 
\begin{lemma}
\label{lemma:subgp}
Given a finite field $\F_q$ and the polynomials \\ 
($\star\star\star$)\hfill   
$x^N-1$ and $c_1+c_2 x^{a_2}+\cdots + c_{t} x^{a_{t}}$, \hfill \mbox{}\\ 
in $\F_q[x]$ with $0\!<\!a_{2}\!<\cdots<\!a_t\!<\!N$,
$\gcd(N,a_2,\cdots,a_t) = 1$, $c_i\!\neq\!0$ for all $i$,  
and $N|(q-1)$, there exists a deterministic $q^{1/4}(\log q)^{O(1)}
+4^t(t\log N)^{O(1)} + t^{\frac{1}{2}+o(1)} N^{\frac{t-2}{t-1}+o(1)} 
(\log q)^{2+o(1)}$ algorithm to decide whether these two polynomials share a 
root in $\F_q$. Furthermore, for some $\delta'|N$ with $\delta'\!\leq\!
\sqrt{t-1}N^{\frac{t-2}{t-1}}$ and $\gamma\!\in\!\{1,\ldots,\delta'\}$, 
the roots of ($\star\star\star$) lie in the union of a set of 
cardinality $2\gamma\sqrt{t-1}N^{\frac{t-2}{t-1}}/\delta'$ and 
the union of $\delta'-\gamma$ cosets of a subgroup of $\F^*_q$ of 
order $N/\delta'$. 
\end{lemma} 

\noindent 
{\bf Proof of Lemma \ref{lemma:subgp}:} 
By Lemma \ref{lemma:expo} we can find an integer $e$ such that, 
if $m_2,\ldots,m_t$ are the unique integers 
in the range $\left[-\lfloor N/2\rfloor,\lceil N/2\rceil\right]$ 
respectively congruent to $ea_2,\ldots,ea_t$, then 
$|m_i|<\sqrt{t-1}N^{\frac{t-2}{t-1}}$ for each $i\!\in\!\{2,\ldots,t\}$.  
Thanks to \cite{mic}, this takes $4^t(t\log N)^{O(1)}$ deterministic 
bit operations. By \cite{shparlinski}, we can then find a generator $\sigma$ 
of $\F^*_q$ within 
$q^{1/4}(\log q)^{O(1)}$ bit operations. For any $\tau\!\in\!\F^*_q$, 
let $\langle \tau\rangle$ denote the multiplicative subgroup of $\F^*_q$ 
generated by $\tau$. 

Now, $x^N-1$ vanishing is the same as 
$x\!\in\!\langle \sigma^{\frac{q-1}{N}}
\rangle$ since $N|(q-1)$. Let $\zeta_{N}\!:=\!\sigma^{\frac{q-1}{N}}$ and 
define $\delta'\!:=\!\gcd(e,N)$. If $\delta'=1$ then the map from $\langle 
\zeta_{N} \rangle$ to 
$\langle \zeta_{N} \rangle$ given by $x\mapsto x^e$ is one-to-one. So finding a 
solution for ($\star\star\star$) is equivalent to finding
$x\!\in\!\langle \zeta_{N} \rangle$ such that
$c_1 +c_2 x^{ea_2}+\cdots + c_{t} x^{ea_{t}}\!=\!0$.   
Thanks to Lemma \ref{lemma:expo}, the last equation can be rewritten as the 
lower degree equation $c_1 +c_2 x^{m_2}+\cdots + c_{t} x^{m_{t}}\!=\!0$, 
and we may conclude our proof by applying Proposition \ref{prop:fac}. 

However, we may have $\delta'\!>\!1$. In which case, the map
from $\langle \zeta_{N}\rangle$ to $\langle\zeta_{N}\rangle$
given by $x\mapsto x^e$ is no longer one-to-one. Instead, 
it sends $\langle\zeta_{N}\rangle$ to a smaller subgroup $\langle
\zeta_{N}^{\delta'}\rangle$ of order $N/\delta'$. We first bound $\delta'$: 
re-ordering monomials if necessary, we may assume that $m_2\!\neq\!0$.  
We then obtain\\ 
\mbox{}\hfill  
$\delta'\!=\!\gcd(e,N)\!\leq\!\gcd(e a_2, N)\!=\!\gcd (m_2,N)\!
\leq\!|m_2|\leq \sqrt{t-1} N^{\frac{t-2}{t-1}}$.\hfill \mbox{} \\   
Any element $x\!\in\!\langle \zeta_{N} \rangle$ can be written
as $\zeta^i_{N} z$ for some $i\!\in\!\{0,\ldots,\delta'-1\}$ and
$z\!\in\!\langle \zeta_{N}^{\delta'}\rangle$. 
It is then clear that $x^N-1\!=\!c_1 +c_2 
x^{a_2}+\cdots +c_{t} x^{a_{t}}\!=\!0$ has a root in $\F^*_q$ if and only if 
there is an  $i\!\in\!\{0,\ldots,\delta'-1\}$ and a 
$z\!\in\!\langle \zeta^{\delta'}_N
\rangle$ with $c_1 +c_2 (\zeta^i_{N}z)^{a_2}
+\cdots + c_{t} (\zeta^i_{N}z)^{a_t} =0$. 
Now, $\gcd(e/\delta',N/\delta')\!=\!1$. 
So the map from $\langle\zeta_{N}^{\delta'} 
\rangle$ to $\langle\zeta_{N}^{\delta'}\rangle$ given by $x\mapsto 
x^{e/\delta'}$ is one-to-one. By the definition of the 
$m_i$, ($\star\star\star$) having a solution is thus equivalent to there 
being an $i\!\in\!\{0,\ldots,\delta'-1\}$ and a 
$z\!\in\!\langle \zeta^{\delta'}_N
\rangle$ with $c_1 +c_2 \zeta^{a_2i}_N z^{m_t/\delta'} 
+\cdots + c_{t} \zeta^{a_ti}_N z^{m_{t}/\delta'} =0$.  
So define the Laurent polynomial  \\ 
\mbox{}\hfill 
$f_i(z)\!:=\!c_1 +c_2 (\zeta^i_{N})^{a_2} z^{m_2/\delta'}
+\cdots + c_{t} (\zeta^i_{N})^{a_{t}} z^{m_{t}/\delta'}$ \hfill \mbox{} \\ 
If $f_i$ is identically zero then
we have found a whole set of solutions for ($\star\star\star$):
the coset $\zeta^i_{N}\langle\zeta_{N}^{\delta'}\rangle$.
If $f_i$ is not identically zero then let
$\ell\!:=\!\min_i { \min(m_i/\delta', 0)}$.
The polynomial $z^{-\ell} f_i(z)$ then has degree bounded from above by
$2\sqrt{t-1} N^{\frac{t-2}{t-1}}/\delta'$. Deciding 
whether the pair of equations
$z^{N/\delta'}-1\!=\!z^{-\ell} f_i(z)\!=\!0$ has a solution for some $i$ 
takes deterministic time 
$\delta'\left(\sqrt{t-1}N^{\frac{t-2}{t-1}}/\delta'\right)^{1+o(1)}
(\log q)^{2+o(1)}$, applying Proposition \ref{prop:fac} $\delta'$ times. 

The final statement characterizing the set of solutions to ($\star\star
\star$) then follows immediately upon defining $\gamma$ to be the 
number of $i\!\in\!\{0,\ldots,\delta'-1\}$ such that $f_i$ is not 
identically zero. In particular, $\gamma\!\geq\!1$ since $\deg f\!<\!N$ and 
thus $f$ is not identically zero on the order $N$ subgroup of $\F^*_q$. \qed
\begin{rem} 
Via fast randomized factoring, we can also pick out 
a representative from each coset of roots within essentially the 
same time bound. Note also that it is possible for some of the  
Laurent polynomials $f_i$ to vanish identically: the polynomial 
$1+x-x^2-x^3$ and the prime $q\!=\!13$, obtained by mimicking Example 
\ref{ex:tetra}, provide one such example (with $\delta'\!=\!6$ and 
$\gamma\!=\!1$). \dia 
\end{rem} 

\smallskip 
We are now ready to prove our first main theorem. \\ 
{\bf Proof of Theorem \ref{thm:meat}:} 
Let $\delta\!:=\!\gcd(q-1,a_2,\ldots,a_{t})$ and 
$y\!=\!x^\delta$. Then the solvability of $f$ is equivalent to
the solvability of the following system of equations: \\ 
\mbox{}\hfill 
$\begin{matrix} c_1 +c_2 y^{a_2/\delta}+\cdots + c_t 
y^{a_t/\delta} =0 \\ y^{\frac{q-1}{\delta}} = 1 \end{matrix}$\hfill \mbox{}\\ 
Since $\gcd\!\left(\frac{a_1}{\delta},\ldots,\frac{a_t}{\delta},
\frac{q-1}{\delta}\right)\!=\!1$, we can solve this problem via Lemma 
\ref{lemma:subgp} (with $N\!=\!\frac{q-1}{\delta}$), within the stated time 
bound. (Note that $q^{1/4}\!\leq\!q^{\frac{t-2}{t-1}}$ for all $t\!\geq\!3$. 
Also, the computation of $\gcd(q-1,a_2,\ldots,a_t)$ 
is dominated by the other steps of the algorithm underlying Lemma 
\ref{lemma:subgp}.) Also, since $y^{\frac{q-1}{\delta}}\!=\!1$, each solution 
$y$ of the preceding $2\times 1$ system induces exactly $\delta$ roots of $f$ 
in $\F_q$. So we can indeed efficiently detect roots of $f$, and the 
second assertion of Lemma \ref{lemma:subgp} gives us the stated 
characterization of the roots of $f$. In particular, $S_2$ is the 
unique order $\frac{q-1}{\delta'}$ subgroup of $\F^*_q$ 
(following the notation of the proof of Lemma \ref{lemma:subgp}). 

The final upper bound then 
follows easily from computing the maximal cardinality of the resulting 
union of cosets, for the cases $\gamma\!\in\!\{1,\eta\}$ (following the 
notation of the proof of Lemma \ref{lemma:subgp}). In particular, 
cosets of $S_2$ do not appear when $\delta'\!=\!1$, and when $\delta'\!>\!1$ 
we clearly have $|S_2|\!\leq\!\frac{q-1}{2}$. \qed 
 
\subsection{The Proof of Corollary \ref{cor:kby1}}  
\label{sub:kby1} 
Deciding whether $0$ is a root of all the $f_i$ is 
trivial, so let us divide all the $f_i$ by a suitable power of $x$ so 
that all the $f_i$ have a nonzero constant term. Next, concatenate 
all the nonzero exponents of the $f_i$ into a single vector of length 
$T\!\leq\!k(t-1)$. Applying Lemma \ref{lemma:expo}, and repeating our 
power substitution trick from our proof of Theorem \ref{thm:meat}, 
we can then reduce to the case where each $f_i$ has degree 
at most $2\sqrt{T}q^{1-T^{-1}}$, at the expense of 
$4^T(T\log q)^{O(1)}$ deterministic bit operations. 

At this stage, we then simply compute 
$g(x)\!:=\!((\cdots(\gcd(\gcd(f_1,f_2),f_3),\ldots),f_k)$ via $k-1$\linebreak 
applications of the Knuth-Sch\"onhage algorithm \cite[Ch.\ 3]{bcs}. This 
takes\\ 
\mbox{}\hfill 
$(k-1)\left(2\sqrt{T}q^{1-T^{-1}}\right)^{1+o(1)}(\log q)^{1+o(1)}$
\hfill\mbox{}\\ 
deterministic bit operations. We then conclude via Proposition \ref{prop:fac}, 
at a cost of\linebreak 
$\left(2\sqrt{T}q^{1-T^{-1}}\right)^{1+o(1)}(\log q)^{2+o(1)}$ big 
operations. 

Summing the complexities of our steps, we arrive at our stated complexity 
bound. \qed 

\section{Hardness in One Variable: Proving Theorems \ref{thm:fphard},  
\ref{thm:adisc}, and \ref{thm:slp}}  
\label{sec:primes} 
\subsection{The Proof of Theorem \ref{thm:fphard}} 
Thanks to Theorem \ref{thm:newplai} we obtain an immediate ZPP-reduction from 
$\sat$ to the detection of roots in $\F_p$ for systems of univariate 
polynomials in $\F_p[x]$. By Lemma \ref{lemma:lin} and Remark 
\ref{rem:lingood} we then obtain a 
BPP-reduction to $2\times 1$ systems. 
Let us now describe a $\zpp$-reduction from $2\times 1$ systems to 
$1\times 1$ systems. 

Suppose $\chi\!\in\!\F_q$ is a
quadratic non-residue. Clearly, the only root in $\F^2_q$ of the
quadratic form $x^2-\chi y^2$ is $(0,0)$. So we can decide the solvability
of $f_1(x)\!=\!f_2(x)\!=\!0$ over $\F_q$ by deciding the 
solvability of $f^2_1-\chi f^2_2$ over $\F_q$. 
Finding a usable $\chi$ is easily done in $\zpp$ via random-sampling
and polynomial-time Jacobi symbol calculation (see, e.g.,
\cite[Cor.\ 5.7.5 \& Thm.\ 5.9.3, pg.\ 110 \& 113]{bs}). 

So there is indeed a 
$\bpp$-reduction from $\sat$ to our main problem, 
and we are done. \qed 

\subsection{The Proof of Theorem \ref{thm:adisc}}  
\label{sub:gcd}   
First note that the hardness of detecting common degree one factors in 
$\F_p[x]$ (or $\bF_p[x]$) for 
{\em pairs} of polynomials in $\F_p[x]$ follows immediately from Theorem 
\ref{thm:newplai} and Lemma \ref{lemma:lin}: the proof of Theorem 
\ref{thm:fphard} above already tells us that there is a $\bpp$-reduction from 
$\sat$ to detecting common roots in $\bF_p$ of pairs of polynomials in 
$\F_p[x]$. Thanks to Assertion (4) of Theorem \ref{thm:newplai}, we also 
obtain a $\bpp$-reduction to detecting common roots, in $\F_p$ instead, 
for pairs of polynomials in $\F_p[x]$. 

So why does this imply hardness for deciding divisibility by the square of a 
degree one polynomial in $\bF_p[x]$ (or $\F_p[x]$)? Assume temporarily 
that Problem (2) is doable in $\bpp$. Consider then, for any 
$f,g\!\in\!\F_p[x]$, the polynomial $H\!:=\!(f+ag)(f+bg)$ where 
$\{a,b\}\!\subset\!\F_p[x]$ is a uniformly random subset of cardinality $2$. 
Note that should $f$ and $g$ have a common factor 
in $\bF_p[x]$, then $H$ has a repeated factor in $\bF_p[x]$. 

On the other hand, if $f$ and $g$ have no common factor, then $f+ag$ and 
$f+bg$ clearly have no common factors. Moreover, thanks to Lemma 
\ref{lemma:sq} and Remark \ref{rem:sqgood}, the 
probability that $f+ag$ and $f+bg$ are both square-free --- and thus 
$H$ is square-free --- is at least 
$\left(1-\frac{2d-1}{q}\right)\left(1-\frac{2d-2}{q}\right)$, assuming $f$ and 
$g$ satisfy the hypothesis of the lemma. 

In other words, to test $f$ and $g$ for common factors, it's enough to  
check square-freeness of $H$ for random $(a,b)$. 

To conclude, thanks to Theorem \ref{thm:newplai}, the pairs of 
polynomials arising from our $\bpp$-reduction from $\sat$ satisfy 
the hypothesis of Lemma \ref{lemma:sq}. Furthermore, thanks to Assertion 
(1) of Theorem \ref{thm:newplai}, our success probability is 
at least $\left(1-\frac{2}{11}\right)^2\!\geq\!\frac{2}{3}$, so we are done. 
\qed  

\subsection{Proving Theorem \ref{thm:slp}}  
We will need the following proposition, due to Ryan Williams.  
\begin{prop} 
\label{prop:will} 
\cite{williams} 
Assume that, for any Boolean circuit of size $L$,  
the Circuit\linebreak Satisfiability Problem can be solved 
in $2^{L-\omega(L)}$ time. Then $\nexp\!\not\subseteq\!\ppoly$. \qed 
\end{prop}

We will also need the following lemma, which is implicit in 
\cite{kipnisshamir}. For completeness, we prove Lemma \ref{lemma:kip} in 
Appendix C. 
\begin{lemma} 
\label{lemma:kip} 
Given a Boolean circuit with $d$ inputs and $L$ gates,
we can find a straight-line program of size $L^{O(1)}$
for a polynomial $f\!\in\!\F_{2^d}[x]$  such
that the circuit is satisfied if and only if
$f$ has a root in $\F_{2^d}$. 
\end{lemma}

\noindent
{\bf Proof of Theorem \ref{thm:slp}:} From Lemma \ref{lemma:kip}, 
an algorithm as hypothesized in Theorem \ref{thm:slp} would imply 
a $2^{L-\omega(L)}$ algorithm for any size $L$ instance of the Circuit 
Satisfiability Problem. By Proposition \ref{prop:will}, we would then 
obtain $\nexp\!\not\subseteq\!\ppoly$. \qed  

\section*{Acknowledgements} 
We would like to thank Igor Shparlinski for insightful comments 
on an earlier draft of this paper. 

\bibliographystyle{amsalpha}

\section*{Appendix A: The Proof of Lemma \ref{lemma:lin}; and Trinomial 
Discriminants}    
Let us first recall the following famous quantitative lemma.
\begin{sz}
Suppose $K$ is any algebraically closed field, 
$f\!\in\!K[x_1,\ldots,x_n]$ is a non-constant 
polynomial of degree $d$,  
and $S\!\subseteq\!K$ has cardinality $N$. Then
$f$ vanishes at no more than $dN^{n-1}$ points of $S^n$. \qed
\end{sz}

\noindent 
{\bf Proof of Lemma \ref{lemma:lin}:} Let $h\!=\!\gcd(f_1,\ldots,f_k)$.
It is then clear that $h\!\in\!\F_q[x]$, $\deg \frac{f_i}{h}\!\leq\!d$
for all $i$, $Z(h)\!=\!Z(f_1,\ldots,f_k)$, and
$Z\!\left(\frac{f_1}{h},\ldots,\frac{f_k}{h}\right)\!=\!\emptyset$.
So if $Z\!\left(\frac{f_1}{h},u_2\frac{f_2}{h}+\cdots+u_k\frac{f_k}{h}\right)
\!=\!\emptyset$ then we clearly obtain
$Z(f_1,u_2f_2+\cdots+u_kf_k)\!=\!Z(f_1,\ldots,f_k)$.
We may thus reduce our lemma to the special case where
$Z(f_1,\ldots,f_k)\!=\!\emptyset$ by simply replacing each $f_i$ by
$\frac{f_i}{h}$. So let us now prove this special case.

Consider the polynomial 
$L(u)\!:=\!\res(f_1,u_2f_2+\cdots+u_kf_k)\!\in\!\F_q[u_2,
\ldots,u_k]$. By construction, for any $\zeta\!\in\!\bF_q$, we either have
$f_1(\zeta)\!\neq\!0$ or $f_i(\zeta)\!\neq\!0$ for some $i\!\geq\!2$. In
the latter case, we see that $u_2f_2(\zeta)+\cdots+u_kf_k(\zeta)\!\neq\!0$
when $u_i\!=\!1$ and all other $u_j$ are $0$. So, by Lemma \ref{lemma:syl},
$L(u)$ is not identically zero. By the Schwartz-Zippel Lemma, we
then obtain that $L(u_2,\ldots,u_k)$ is nonzero for at least a fraction of
$1-\frac{d}{q}$ of the $(u_2,\ldots,u_k)\!\in\!\F^{k-1}_q$.
Moreover, Lemma \ref{lemma:syl} tells us that at any such point,
$Z(f_1,u_2f_2+\cdots+u_kf_k)\!=\!\emptyset$. So we are done. \qed 

We now make a final observation about the roots of trinomials over 
finite fields, easily following from \cite[Lemma 5.3]{airr}. 
\begin{lemma} 
\label{lemma:tri} 
Suppose $f(x)\!=\!c_1+c_2x^{a_2}+c_3x^{a_3}\!\in\!\F_q[x]$ has 
degree $<\!q$, $\cA\!:=\!\{0,a_2,a_3\}$,\linebreak $0\!<\!a_2\!<\!a_3$, and 
$\gcd(a_2,a_3)\!=\!1$. Recall that $\zeta\!\in\!\bF_q$ is a degenerate root 
of $f \Longleftrightarrow$\linebreak 
$f(\zeta)\!=\!f'(\zeta)\!=\!0$. Then: 

\smallskip 
\noindent
(0) $\Delta_\cA(f)=(a_3-a_2)^{a_3-a_2}a^{a_2}_2c^{a_3}_2-
(-a_3)^{a_3}c^{a_3-a_2}_1c^{a_2}_3$.\\ 
(1) $\Delta_\cA(f)\!\neq\!0\Longleftrightarrow f$ has no degenerate
roots in $\bF_q$. In which case, we also have\\  
\mbox{}\hspace{.75cm}\scalebox{.95}[1]
{$\Delta_\cA(f)\!=\frac{(-1)^{a_3}c^{a_2-1}_3}{c^{a_2-1}_1}
\!\!\!\!\prod\limits_{f(\zeta)=0}f'(\zeta)$ where the product 
ranges over the $a_3$ distinct roots of $f$ in $\bF_q$.}\\ 
(2) Deciding whether $f$ has a degenerate root in $\bF_p$ can be  
done in time polynomial in\linebreak 
\mbox{}\hspace{.75cm}$\log q$.\\  
(3) If $f$ has a degenerate root $\zeta\!\in\!\bF^*_p$ then  
$(\zeta^{a_2},\zeta^{a_3})\!=\!\frac{c_1}{a_3-a_2}
\left(-\frac{a_3}{c_2},\frac{a_2}{c_3}\right)$. 
In particular, such a\linebreak 
\end{lemma} 

\section*{Appendix B: The Proof of Lemma \ref{lemma:sq}}    
For $2d-1\!\geq\!p$ the lemma is vacuous, so let us assume $2d-1\!<\!p$. 
Note also that the polynomial $f+ag$ is irreducible in $\F_p[x,a]$, 
since $f$ and $g$ have no common factors in $\F_p[x]$. The splitting 
field $L\!\subsetneqq\!\overline{\F_p(a)}$ of $f(x)+ag(x)$ must 
have degree $[L:\F_p(a)]$ dividing $(\deg f)!$. Since $\deg f\!\leq\!d\!<\!p$, 
$p$ can not divide $[L:\F_p(a)]$ and thus $L$ is a separable extension of 
$\F_p(a)$, i.e., $f+ag$ has no degenerate roots in $\overline{\F_p(a)}$. So 
the classical discriminant of $f+ag$ (where the coefficients are considered 
as polynomials in $a$) is a polynomial in $a$ that is not identically 
zero. Furthermore, from Definition \ref{dfn:syl}, 
$\res_{(d,d-1)}(f+ag,f'+ag')\!\in\!\F_p[a]$ 
has degree at most $d+d-1\!=\!2d-1$. So by Lemma \ref{lemma:lin}, 
the classical discriminant of $f+ag$ is non-zero for at least 
$1-\frac{2d-1}{p}$ of the $a\!\in\!\F_p$. Thanks to Lemma \ref{lemma:syl}, 
we thus obtain that $f+ag$ is square-free for at least a fraction of 
$1-\frac{2d-1}{p}$ of the $a\!\in\!\F_p$. \qed 

\section*{Appendix C: The Proof of Lemma \ref{lemma:kip}}   
A Boolean circuit can be viewed as a straight-line program using
Boolean variables and Boolean operations. One can replace the Boolean 
operations by polynomials over $\F_2$: \\ 
\mbox{}\hfill $x_1 \wedge x_2 =  x_1x_2$\hfill \mbox{}\\ 
\mbox{}\hfill $x_1 \vee x_2 = x_1 + x_2 + x_1 x_2$\hfill \mbox{}\\
\mbox{}\hfill $\neg x_1 = 1 - x_1$\hfill \mbox{}\\
Hence a straight-line program for a Boolean function of size $L$ with 
$d$ inputs can be converted into a straight-line program for a polynomial
$f(x_0, x_1, \cdots, x_{d-1}) \in \F_2[x_0, x_1, \cdots, x_{d-1}]$
of size $O(L)$.

Let $b(x)$ be an irreducible polynomial of degree $d$ over $\F_2$.
Let $\alpha$ be one root of $b(x)$. Then
$\{1,\alpha,\alpha^2,\ldots,\alpha^{d-1}\}$ is a basis 
for $\F_{2^d}$ over $\F_2$. Then any element $x\!\in\!\F_{2^d}$ can be 
written uniquely as $x\!=\!x_0+x_1\alpha +\cdots +x_{d-1}\alpha^{d-1}$, 
where $x_i \in \F_2$ for all $i$. So we obtain the system of linear 
equations\\ 
\mbox{}\hfill 
$\begin{bmatrix}
1 & \alpha & \cdots & \alpha^{d-1}\\
1 & \alpha^2 & \cdots & \alpha^{2(d-1)}\\
1 & \alpha^4 & \cdots & \alpha^{4(d-1)}\\
\vdots\\
1 & \alpha^{2^{d-1}} & \cdots & \alpha^{2^{d-1}(d-1)}\\
\end{bmatrix}\begin{bmatrix}x_0\\ x_1\\ x_2\\ \vdots \\ x_{d-1}\end{bmatrix} 
=\begin{bmatrix}x\\ x^2\\ x^4\\ \vdots \\ x^{2^{d-1}}\end{bmatrix}$.\hfill 
\mbox{}\\ 
The underlying matrix is Vandermonde and thus non-singular.
So we can represent each $x_i$ as a linear combination of
$x,x^{2^1},x^{2^2},\ldots,x^{2^{d-1}}$ over $\F_{2^d}$.
Replacing each $x_i$ by the appropriate linear combination of 
high powers of $x$, in the SLP for $f$, we obtain our lemma. \qed 

\end{document}